\documentclass{amsart}

\usepackage{enumerate}
\theoremstyle{plain}
\newtheorem{Thmf}{Th\'eor\`eme}
\newtheorem{Propf}[Thmf]{Proposition}
\newtheorem{Thm}{Theorem}[section]

\newtheorem{Prop}[Thm]{Proposition}
\newtheorem{Lem}[Thm]{Lemma}

\numberwithin{equation}{section}
\DeclareMathOperator{\M}{M}

\DeclareMathOperator{\tr}{tr}
\DeclareMathOperator{\im}{im}

\DeclareMathOperator{\dom}{dom}

\begin{document}

\title[A Rationality Criterion]
{A Rationality Criterion For Unbounded Operators}

\author[P. A. Linnell]{Peter A. Linnell}

\address{Math \\ VPI \\ Blacksburg \\ VA 24061--0123
\\ USA}

\email{linnell@math.vt.edu}
\urladdr{http://www.math.vt.edu/people/linnell/}

\begin{abstract}
Let $G$ be a group, let $U(G)$ denote the set of unbounded operators
on $L^2(G)$ which are affiliated to the group von Neumann algebra
$W(G)$ of $G$, and let $D(G)$ denote the division closure of $\mathbb
{C}G$ in $U(G)$.  Thus $D(G)$ is the smallest subring of $U(G)$
containing $\mathbb {C}G$ which is closed under taking inverses.
If $G$ is a free group then $D(G)$ is a division ring, and in this
case we shall give a criterion for an element
of $U(G)$ to be in $D(G)$.  This extends a
result of Duchamp and Reutenauer, which was concerned with proving a
conjecture of Connes.

\begin{center}
\textbf{Un crit\`ere de rationalit\'e pour op\'erateurs non born\'es}
\end{center}

\textbf{R\'esum\'e } --
Soient $G$ un groupe, $U(G)$ l'ensemble d'op\'erateurs non born\'es
affili\'es \`a l'alg\`ebre de von Neumann de groupe de $G$,
et $D(G)$ la
cl\^oture de division de $\mathbb {C}G$ dans $U(G)$.
Ainsi $D(G)$ est le plus petit anneau qui est ferm\'e
sous l'op\'eration d'inverse.  Si $G$
est un group libre, nous donnons un crit\`ere pour qu'un
\'el\'ement de $U(G)$ soit dans $D(G)$.
\end{abstract}

\keywords{unbounded operator, finite rank, division closure}

\subjclass{Primary: 22D25; Secondary: 20C07, 46L10}

\maketitle

\textbf{Version fran\c{c}ais abr\'eg\'ee } --
Soient
$G$ un groupe libre, $L^2(G)$ l'espace de Hilbert avec base
orthonorm\'ee $\{g \mid g \in G\}$,
$C_r^*(G)$ l'alg\`ebre r\'eduite de groupe de $G$,
$W(G)$ l'alg\`ebre de von Neumann de
groupe de $G$, $U(G)$ l'ensemble d'op\'erateurs ferm\'es,
affili\'es \`a $W(G)$, $D(G)$ la
cl\^oture de division de $\mathbb {C}G$ dans $U(G)$,
et $S(G)$ la
cl\^oture de division de $\mathbb {C}G$ dans $C_r^*(G)$.
Ainsi $D(G)$ est le plus petit anneau qui est ferm\'e
sous l'op\'eration d'inverse.
J'ai montr\'e dans \cite{div} que $D(G)$ est un anneau de division.
Il y a un $G$-ensemble libre et un
op\'erateur unitaire $P \colon L^2(G) \to L^2(E) \oplus
\mathbb {C}$ tel que l'op\'erateur $P \alpha - \alpha P$
soit de rang fini pour tout $\alpha \in \mathbb {C}G$
\cite[p.~341]{connes}.  Si $A$ est une sous-alg\'ebre de $U(G)$,
appelons
$A_{\text{fin}}$ la plus grande sous-alg\'ebre de $A$ telle
que $P \alpha -\alpha P$ soit de rang finie pour tout $\alpha \in A$.
Reutenauer et Duchamp ont montr\'e dans \cite{duchreut} que la
cl\^oture de division de $\mathbb {C}G$ dans $C_r^*(G)$ est
$(W(G))_{\text{fin}}$;
ceci a r\'epondu \`a une question de Connes \cite[p.~342]{connes}.
Nous \'etendrons ce r\'esultat au $U(G)$.

D\'efinissons les sous-ensembles $R(G)$ et $R'(G)$
de $G$ comme suit.  Pour $u \in U(G)$,
nous disons que $u \in R(G)$ si et seulement si toutes les fois que $u
= s^{-1}a = bt^{-1}$ avec $a,b,s,t \in W(G)$, alors $sPb - aPt$ et
$sP^{-1}b - tP^{-1}a$ soit de rang fini, alors que
nous disons que $u \in R'(G)$ si et seulement si nous pouvons
\'ecrire $u
= s^{-1}a = bt^{-1}$ avec $a,b,s,t \in W(G)$
et tels que $sPb - aPt$ et
$sP^{-1}b - tP^{-1}a$ ont rang fini.
Nous pouvons maintenant \'enoncer
\begin{Thmf} \label{Tf}
$D(G) = R(G) = R'(G)$ et $D(G) \cap W(G) = S(G)$.  En outre si $u \in
D(G)$, alors nous pouvons \'ecrire $u = s^{-1}a = bt^{-1}$ avec
$a,b,s,t \in S(G)$.
\end{Thmf}
Ainsi en particulier, chaque \'el\'ement de $D(G)$ peut \^etre
\'ecrit sous la forme $s^{-1}a$ avec $a,s \in C_r^*(G)$.
La d\'emonstration du Th\'eor\`eme \ref{Tf} d\'epend crucialement des
r\'esultats de \cite{duchreut}.
Une autre description de $D(G)$ est donn\'ee par le r\'esultat
suivant.
\begin{Propf}
Soient $u \in U(G)$.  Alors $u \in D(G)$ si et seulement s'il y a
un sous-espace $M$ de codimension finie dans $L^2(G)$
tels que les restrictions de $Pu$ et $uP$ \`a $M$  sont \'egales.
\end{Propf}

\section{Introduction} \label{Sintroduction}

Let $H$ be a Hilbert space and for
$u,v \in H$, let
$\langle u,v \rangle$ indicate the inner product of $u$ and $v$.
The set of all closed densely defined linear operators acting on the
left of $H$ will be denoted by $\mathcal {U}(H)$, and the subset
consisting of bounded operators will be denoted by $\mathcal {B}(H)$.
The adjoint $\theta^*$ of $\theta \in \mathcal {U}(H)$ satisfies
$\langle \theta u, v \rangle = \langle u, \theta^*v \rangle$
whenever $\theta u$ and $\theta^*v$ are defined.  Now
let $G$ be a group and let $L^2(G)$ denote the Hilbert space with
Hilbert basis $\{ g \mid g \in G\}$.  Thus $L^2(G)$ consists of all
formal sums $\sum_{g\in G}a_g g$ where $a_g \in \mathbb {C}$ and
$\sum_{g \in G} |a_g|^2 < \infty$, and has inner product defined by
\[
\langle \sum_{g \in G} a_g g, \sum_{h \in G} b_h h \rangle = \sum_{g
\in G} a_g \bar{b}_g
\]
where $\bar{\ }$ denotes complex conjugation.  If $\alpha = \sum_{g \in
G} a_g g \in \mathbb {C}G$ (so $a_g \in \mathbb {C}$ and $a_g = 0$
for all but finitely many $g$) and $\beta = \sum_{g\in G} b_gg
\in L^2(G)$, then
\[
\alpha \beta = \sum_{g,h \in G} a_g b_h gh = \sum_{g\in G} (\sum_{h
\in G} a_{gh^{-1}} b_h ) g \in L^2(G)
\]
and the map $\beta \mapsto \alpha \beta $ (left multiplication by
$\alpha$) is in $\mathcal {B}(L^2(G))$.  It follows that we may
identify $\mathbb {C}G$ with a subring of $\mathcal {B}(L^2(G))$.
By definition the reduced group C*-algebra $C_r^*(G)$ of $G$ is the
norm closure of $\mathbb {C}G$ in $\mathcal {B}(L^2(G))$, and
the group von Neumann algebra $W(G)$ of $G$ is the weak
closure of $\mathbb {C}G$ in $\mathcal {B} (L^2(G))$; thus
$\mathbb {C}G \subseteq C_r^*(G) \subseteq W(G)$ and $W(G)$ is a
finite von Neumann algebra.  Let $U(G)$ denote the operators
in $\mathcal {U}(G)$ which are affiliated to $W(G)$
\cite[p.~150]{ber}.  Then $U(G) = U(G)^*$, $U(G)$ is a $*$-regular
ring \cite[definition 1 of \S 51]{berbook}
containing $W(G)$, and every element of $U(G)$ can be written in
the form $s^{-1}a$ and also $as^{-1}$, where $a \in W(G)$ and $s$ is
a nonzero divisor in $W(G)$ \cite[theorems 1 and 10]{ber}.
Using the fact that $U(G)$ is a
$*$-regular ring, we see that $s$ is a nonzero divisor in $W(G)$ if
and only if $sw \ne 0$ (or $ws \ne 0$) whenever $0 \ne w \in W(G)$.
Furthermore any finite set of elements
in $U(G)$ have a common denominator, so for example if $u,v \in U(G)$,
then there exist $a,b,s \in W(G)$ with $s$ a nonzero divisor such that
$u = s^{-1}a$ and $v = s^{-1} b$.

Suppose $H$ is the Hilbert space sum of an arbitrary number (finite or
infinite) of copies of $L^2(G)$.  For $h \in H$, we shall write $h =
(h_1, h_2, \dots )$ where $h_i$ denotes the $i$th component of $h$.
If $\theta \in \mathcal {U}(L^2(G))$,
then $\theta$ defines an element of $\mathcal {U}
(H \oplus \mathbb {C})$
according to the rule $\theta (h_1, h_2, \dots, c) = (\theta h_1,
\theta h_2, \dots, c)$ where $c \in \mathbb {C}$.  The domain of this
operator is $\dom (\theta) \times \dom (\theta) \times \dots \times
\mathbb {C}$.  If $\theta \in U(G)$ and $\theta = s^{-1}a = bt^{-1}$
where $a,b,s,t \in W(G)$ and $s,t$ are nonzero divisors,
then $tL^2(G) \subseteq
\dom (\theta) = \{h \in L^2(G) \mid ah \in sL^2(G) \}$.

Suppose $R$ is a subring of the ring $S$.  Then the division closure
of $R$ in $S$ is the smallest subring $D$ of $S$ containing $R$ which
is closed under taking inverses in $S$ (i.e.\ $d \in D$ and $d$
invertible in $S$ implies $d^{-1} \in D$).  We shall let $D(G)$ denote
the division closure of $\mathbb {C}G$ in $U(G)$.
Clearly if there is a
division ring $E$ such that $R \subseteq E \subseteq S$, then the
division closure of $R$ in $S$ is a division subring of $E$.
Also the rational closure of $R$ in $S$
is the subset $T$ of $S$ defined by the property $t \in T$ if
and only if there exists an integer $n$ and $M \in \M_n(R)$ such
that $M$ is invertible in $\M_n(S)$ and $t$ is one of the entries of
$M^{-1}$ \cite[p.~382]{cohn}.  The rational closure of $R$ in $S$ is
always a subring of $S$ \cite[theorem 7.1.2]{cohn} which contains the
division closure of $R$ in $S$ \cite[exercise 7.1.4]{cohn}.  Often
the division closure is equal to the rational closure, but there are
examples when the division closure is strictly contained in the
rational closure.  However in the case that the division closure is a
division ring, then it is clear that the division closure is equal to
the rational closure.

For the rest of this paper $G$ will be a free group.
By \cite[theorem 1.3]{div}, $D(G)$ is a division ring and so $D(G)$ is
equal to the rational closure of $\mathbb {C}G$ in $U(G)$.
If $G$ is free on the set $X$, then it is shown on page 573 of
\cite{div} that $D(G)$ is isomorphic to the free
field on $X$ over $\mathbb {C}$ \cite[p.~224]{cohn2}.

Let $G$ act on the one element set $\{*\}$ according to the rule
$g* = 0$ for all $g \in G$.  Then there is a free left $G$-set $E$ and a
bijection $\pi \colon G \to E \cup \{*\}$ such that $\{b \in G \mid
\pi gb \ne g \pi b \}$ is finite for all $g \in G$
\cite[ p.~341]{connes} ($E$ here is $T^1$ there).
Then $\pi$ extends to a unitary operator $P \colon L^2(G)
\to L^2(E) \oplus \mathbb {C}$ with the
property that $P a - a P$ has finite rank (i.e.\ $\im (Pa - aP
)$ has finite dimension over $\mathbb {C}$) for all $a \in \mathbb
{C}G$; this can be seen from the proof of
\cite[lemma IV.5.1(a) on p.~342]{connes},
where $P$ there is the same as $P$
here.  For any subalgebra $A$ of $\mathcal {B}(L^2(G))$, let
$A_{\text{fin}} = \{a \in A \mid Pa - aP$ has finite rank\} so if
$A \supseteq \mathbb {C}G$, then $A_{\text{fin}} \supseteq \mathbb
{C}G$.  Let $S(G)$ denote the rational closure of $\mathbb {C}G$ in
$C_r^*(G)$.
Then \cite[remark 3 on p.~342]{connes} shows that $S(G) \subseteq
(C_r^*(G))_{\text{fin}}$, and the question of whether $S(G) =
(C_r^*(G))_{\text{fin}}$ is posed there.
This question was answered in the affirmative by
\cite[th\'eor\`eme 7]{duchreut}, where the stronger result, that the
division closure of $\mathbb {C}G$ in $C_r^*(G)$ is equal to
$(W(G))_{\text{fin}}$, was proved so in particular $S(G)$ is also the
division closure of $\mathbb {C}G$ in $C_r^*(G)$.
Christophe Reutenauer has told
me that when he and Duchamp proved this, Connes posed the
problem of extending their result to $U(G)$.
The purpose of this paper is to give an answer to this problem, and
then to give a few simple applications of the result.

Define subsets $R(G), R'(G)$ of $U(G)$ as follows.  For $u \in
U(G)$, we say that $u \in R(G)$ if and only if
whenever $u = s^{-1}a = bt^{-1}$ with
$a,b,s,t \in W(G)$, then $s P b - a P t$ and $s P ^{-1}b -a
P^{-1} t$ have finite rank (the former is a bounded linear operator
$L^2(G) \to L^2(E) \oplus \mathbb {C}$, and the latter is a bounded
linear operator $L^2(E) \oplus \mathbb {C} \to L^2(G)$), while we 
say that $u \in R'(G)$ if and only if we may write $u = s^{-1}a = 
bt^{-1}$ with $a,b,s,t \in W(G)$ and such that $s P b - a P t$
and $s P ^{-1}b -a P^{-1} t$ have finite rank.  As remarked above,
there is always at least one way to write $u = s^{-1}a = bt^{-1}$ 
with $a,b,s,t \in W(G)$, consequently $R(G) \subseteq R'(G)$.
We can now state
\begin{Thm} \label{Tmain}
$D(G) = R(G) = R'(G)$ and $D(G) \cap W(G) = S(G)$.
Furthermore if $u \in D(G)$, then we may write
$u = s^{-1}a = bt^{-1}$ with $a,b,s,t \in S(G)$.
\end{Thm}
It is easy to read off a number of consequences of this result, for
example we can now state that
every element of $D(G)$ can be written in the form $s^{-1}a$ with
$a,s \in (C_r^*(G))_{\text{fin}}$.  It seems plausible that the
definition of $R(G)$ could be weakened to requiring only that $sP b
- aP t$ has finite rank, but I have been unable prove this.
Theorem \ref{Tmain} generalizes
\cite[th\'eor\`eme 7]{duchreut}, and the proof depends crucially on
the results of \cite{duchreut}.

Finally we give two other ways of defining $R(G)$.  For the first
let
\[
F =
\begin{pmatrix}
0&P^{-1}\\
P&0
\end{pmatrix}
\in \mathcal {B}(L^2(G) \oplus L^2(E) \oplus \mathbb {C}),
\]
so $F$ yields a Fredholm module as described on page 341 of
\cite{connes}.  Then for $u \in U(G)$, we can say that $u \in R(G)$ if
and only if whenever $u = s^{-1}a = bt^{-1}$ with $a,b,s,t \in W(G)$,
then $s F b - aFt$ has finite rank; this is obvious.
The second way is described by the following result.
\begin{Prop} \label{Pmain}
Let $u \in U(G)$.  Then $u \in D(G)$ if and only if there exists a
subspace $M$ of finite codimension in $L^2(G)$ such that the
restrictions of $Pu$ and $uP$ to $M$ are equal.
\end{Prop}
Of course $Pu$ and $uP$ will \emph{not} in general be bounded
operators; when we say that two unbounded operators are equal, then we
implicitly assume that there domains of definition are equal.

I am very grateful to Christophe Reutenauer for introducing me to the
problem studied in this paper, and for some useful discussions.

\section{Notation, Terminology and Assumed Results} \label{Snotation}

Most of the notation and terminology used in this paper has already
been defined above.  Mappings will be written on the left and $\mathbb
{C}$ will denote the complex numbers.  All rings will have a 1,
and subrings will have the same 1.  If $E$ is a set, then $L^2(E)$
will denote the Hilbert space with Hilbert basis $\{e \mid e \in
E\}$.  We shall let $\im \theta$ and $\ker \theta$ denote
the image and kernel of the map $\theta$ respectively.
If $n$ is a positive integer, then $\M_n(R)$ will indicate the $n$
by $n$ matrices over a ring $R$.
A projection in $\mathcal {B}(H)$ is
an element $e$ such that $e = e^* = e^2$.  If $u$ is an unbounded
operator, then $\dom (u)$ will indicate the domain of $u$, in other
words the subspace on which it is defined.
We need the following three elementary lemmas.

\begin{Lem} \label{Linvariance}
Let $\theta \colon H \to K$ and $\phi \colon K \to L$ be bounded
linear maps between Hilbert spaces.
\begin{enumerate}[\normalfont (i)]
\item If $\ker \phi = 0$ and $\phi \theta $ has finite rank, then
$\theta $
also has finite rank.

\item If $\im \theta $ is dense in $K$ and $\phi \theta $ has finite
rank,
then $ \phi $ also has finite rank.
\end{enumerate}
\end{Lem}
\begin{proof}
(i) is obvious.  For (ii), since $\theta H$ is dense in $K$, we see
that
$\phi \theta H$ is dense in $\phi K$.  But $\phi \theta  H$ is finite
dimensional and therefore closed,
hence $\phi K$ is finite dimensional and the result is proven.
\end{proof}
\begin{Lem} \label{Ldense}
Let $\theta \in W(G)$.  If $\theta$ is a nonzero divisor, then
$\ker \theta = 0$ and $\im \theta$ is dense in $L^2(G)$.
\end{Lem}
\begin{proof}
Since $\theta$ is a nonzero divisor in $W(G)$, it is invertible in
$U(G)$ and it follows that $\ker \theta = 0$.  Also $\theta^*$ is a
nonzero divisor in $W(G)$, and we deduce from this that $\im \theta$
is dense in $L^2(G)$.
\end{proof}
\begin{Lem} \label{Lu-lambda}
Let $u \in U(G)$.  Then there exists $\lambda \in \mathbb {C}$ such
that $u - \lambda$ is invertible in $U(G)$.
\end{Lem}
\begin{proof}
For each $\lambda \in \mathbb {C}$, let $K_{\lambda} = \{x \in W(G)
\mid (u-\lambda)x = 0\}$, a right ideal of $W(G)$.  We first show that
$K_{\lambda} = 0$ for some $\lambda \in \mathbb {C}$.  Since $W(G)$ is
a von Neumann algebra, there is a unique projection $e_{\lambda} \in
W(G)$ such that $e_{\lambda}W(G) = K_{\lambda}$.  Let $\tr \colon W(G)
\to \mathbb {C}$ denote the trace map, as described for example in
\cite[p.~352]{zero}.  If we write an element $\alpha$ of $W(G)$ in the
form $\sum_{g \in G} a_g g$ where $a_g \in \mathbb {C}$, then $\tr
\alpha = a_1$.  Also if $e$ is a nonzero projection in $W(G)$,
then $0 < \tr
e \le 1$.  Note that the sum $\sum_{\lambda} K_{\lambda}$ is direct,
so by using \cite[lemma 12]{zero}, we see that $\sum_{\lambda \in
S} \tr(e_{\lambda}) \le 1$ for any finite subset $S$ of $\mathbb
{C}$.  It follows that the number of $\lambda$ for which $e_{\lambda}
\ne 0$ is countable, and we deduce that there exists $\lambda \in
\mathbb {C}$ (in fact uncountably many such $\lambda$) such that
$e_{\lambda} = 0$.  For this $\lambda$, we have $K_{\lambda} = 0$ and
$u - \lambda$ is a nonzero divisor in $W(G)$.
Since every element of $U(G)$ can be written in the form $s^{-1}a$
and also $as^{-1}$
with $a,s \in W(G)$, we conclude that $u - \lambda$ is a nonzero
divisor in $U(G)$.  But every element of $U(G)$ is either a zero
divisor or invertible, and the result follows.
\end{proof}

\section{Proofs}

It will be clear from the next lemma that
$\mathbb {C}G \subseteq R(G) = R'(G)$.
We are going to show that
$R(G)$ is a subring which is closed under taking inverses and
adjoints.  This will mean in particular that $R(G)$ is division closed
and so will contain the division closure of $\mathbb {C}G$ in $U(G)$.
First we
show that we need only check the condition $sP b - a P t$,
$sP^{-1} b - a P^{-1} t$ have finite rank
for one choice of $a, b, s, t$ satisfying $u =
s^{-1}a = bt^{-1}$.

\begin{Lem} \label{Ldefined}
Let $u \in U(G)$ and suppose $u =  s^{-1}a = bt^{-1}$, where $a,b,s,t
\in W(G)$ and $s,t$ are nonzero divisors.  If $s P b - a P t $
and $s P^{-1} b - a P^{-1} t $ have finite rank,
then $u \in R(G)$.
\end{Lem}
\begin{proof}
Suppose $u = s_1^{-1} a_1$.  Then we need to show that $s_1 P b
- a_1 P t$ and $s_1 P ^{-1} b - a_1 P ^{-1} t $ have finite
rank.
There are nonzero divisors $x,x_1 \in W(G)$ such that $s s_1^{-1} =
x^{-1} x_1$.  Then $xs = x_1 s_1$ and $xa = x_1 a_1$, hence
\[
x_1 (s_1 P b - a_1 P t) = x (s P b -a P t)
\]
and we see that $x_1 (s_1 P b - a_1 P t)$ has finite rank.
Using Lemmas \ref{Linvariance} and \ref{Ldense},
we deduce that $s_1 P b - a_1P t$
has finite rank.  Similarly $s_1 P^{-1} b - a_1 P ^{-1}
t $ has finite rank.
If $u = b_1 t_1^{-1}$, then in a similar
fashion we can show that $s_1 P b_1 - a_1 P t_1$
and $s_1 P^{-1} b_1 - a_1 P ^{-1} t_1$ have finite rank.
This establishes the result.
\end{proof}
Next we show that $R(G)$ is closed under addition.
\begin{Lem} \label{Ladd}
Let $u, v \in R(G)$.  Then $u + v \in  R(G)$.
\end{Lem}
\begin{proof}
Write $u = s^{-1}a = bt^{-1}$
and $v = s^{-1}c = dt^{-1}$, where $a,b,c,d,s,t \in W(G)$ and $s,t$
are nonzero divisors.  Then
$s P b - a P t $,
$s P^{-1} b - a P^{-1} t $,
$s P d - c P t $,
$s P^{-1} d - c P^{-1} t$ have finite rank and
$u + v = s^{-1} (a + c) = (b + d) t^{-1}$,
consequently
\[
s P (b+d) - (a + c) P t = (s P b - a P t) + (s P d - c P
t)
\]
has finite rank.  Similarly
\[
s P^{-1} (b+d) - (a + c) P^{-1} t =
(s P^{-1} b - a P^{-1} t) + (s P^{-1} d - c P^{-1} t)
\]
has finite rank.
Using Lemma~\ref{Ldefined}, we deduce that $u + v \in R(G)$.
\end{proof}
Now we show that $R(G)$ is closed under multiplication.
\begin{Lem} \label{Lmult}
Let $u,v \in R(G)$.  Then $uv \in R(G)$.
\end{Lem}
\begin{proof}
By Lemma~\ref{Lu-lambda}, there exist $\lambda, \mu \in \mathbb {C}$ such
that $u - \lambda$, $v - \mu $ are invertible in $U(G)$, so using
Lemma~\ref{Ladd}, we may assume that $u,v$ are invertible in $U(G)$.
This means when we write $u = s^{-1} a$ with $a,s \in W(G)$, not only
$s$ but also $a$ are nonzero divisors in $W(G)$.  Write $v = t^{-1}b$
where $b, t$ are nonzero divisors in $W(G)$, and then write $a t^{-1}
= w^{-1} c$ where $c,w$ are nonzero divisors in $W(G)$.  Then
\[
uv = s^{-1}a t^{-1}b = (ws)^{-1}(wa)(ct)^{-1}(cb)
\]
and $wa = ct$.  Thus we may write $u = p^{-1}q$, $v = q^{-1} r$ where
$p,q,r$ are nonzero divisors in $W(G)$, and similarly we may write $u
=xy^{-1}$, $v = yz^{-1}$ where $x,y,z$ are nonzero divisors in $W(G)$.
Then $uv = p^{-1} r = x z^{-1}$.
Since $u,v \in R(G)$, we have
\[
p P x - q P y, \ q P y - r P z
\]
have finite rank and hence $p P x - r P z =
(p P x - q P y) - ( q P y - r P z)$
has finite rank.  Similarly
$p P^{-1} x - r P^{-1} z $ has finite rank
and an application of Lemma \ref{Ldefined} completes the proof.
\end{proof}
Now we show that $R(G)$ is closed under taking inverses.
\begin{Lem} \label{Linv}
Let $u \in R(G)$.  If $u$ is invertible in $U(G)$, then $u^{-1} \in
R(G)$.
\end{Lem}
\begin{proof}
Write $u = s^{-1}a = bt^{-1}$ where $a,b,s,t \in W(G)$, all nonzero
divisors because $u$ is invertible in $U(G)$.   Then $u^{-1} =
a^{-1}s = tb^{-1}$.  Since $u \in R(G)$, we know that
\[
s P b - a P t, \ s P^{-1} b - a P ^{-1} t
\]
have finite rank.  Therefore
\[
a P t  - s P b, \ a P^{-1} t - s P^{-1} b
\]
have finite rank.  The result now follows from Lemma \ref{Ldefined}.
\end{proof}
Finally we show that $R(G)$ is closed under the adjoint operation.
\begin{Lem} \label{Ladj}
Let $u \in R(G)$.  Then $u^* \in R(G)$.
\end{Lem}
\begin{proof}
Write $u = s^{-1} a = bt^{-1}$, where $a,b,s,t \in W(G)$ and $s,t$ are
nonzero divisors.  Then $u^* = (t^*)^{-1} b^* = a^* (s^*)^{-1}$.
Since $u \in R(G)$, we know that $s P b - a P t$ and $s P^{-1} b
- a P ^{-1} t$ have finite rank.  If $T$ is a bounded linear
map between Hilbert spaces
with finite rank, then $T^*$ also has finite rank.
Furthermore $P^* = P^{-1}$ because $P$ is a unitary operator.
Therefore $t^* P^{-1} a^* - b^*
P^{-1} s^* $ and $t^* P a^* - b^* P s^*$ have finite rank.  The
result now follows from Lemma \ref{Ldefined}.
\end{proof}
It now follows from Lemmas \ref{Ldefined}, \ref{Ladd}, \ref{Lmult},
\ref{Linv} and \ref{Ladj} that $R(G)$ is a subring of $U(G)$ containing
$D(G)$ which is closed under the * operation and taking inverses.

\begin{proof}[Proof of Theorem \ref{Tmain}]
We have already shown that $R(G) = R'(G)$, and we
see from \cite[propositions 5 and 9]{duchreut}
and Lemma \ref{Ldefined} that $R(G) \cap W(G)
= S(G)$.  Now let $u \in R(G)$ and set
$a = (1 + uu^*)^{-1} u$ and $s = (1 + uu^*)^{-1}$.  Then $a,s$ are
well defined elements of $W(G)$ and $s$ is a nonzero divisor, by
\cite[15.12.6]{dieu2}.  Since $R(G)$ is a subring of $U(G)$ closed
under taking inverses and adjoints, we see that $a,s \in R(G) \cap
W(G)$.  The result follows.
\end{proof}

\begin{proof}[Proof of Proposition \ref{Pmain}]
First suppose $u \in D(G)$.  Then by Theorem \ref{Tmain}, there exist
$a,s \in S(G)$ such that $u = s^{-1}a$.
Then $\dom (Pu) =
\{x \in L^2(G) \mid ax \in sL^2(G)
\}$ and $\dom (uP)= \{ x \in L^2(G) \mid a P x \in
s PL^2(G) \}$.  Since $Pa -aP$ and $Ps - sP$ have finite rank,
there are subspaces $M_1$ and $M_2$ of finite codimension in $L^2(G)$
such that $Pa -a P$ is zero on $M_1$ and
$Ps - sP$ is zero on $M_2$.
Then $P^{-1}sP M_2 = sM_2$, hence there are subspaces $N_1, N_2$ of
finite codimension in $L^2(G)$ such that $N_1 \cap sL^2(G) \subseteq
sM_2$ and $N_2 \cap P^{-1}sPL^2(G) \subseteq sL^2(G)$.  Now choose
subspaces $M_3, M_4$ of finite codimension in $L^2(G)$ such that
$aM_3 \subseteq N_1$ and $aM_4 \subseteq N_2$, and set $M = M_1 \cap
M_3 \cap M_4$.

Suppose $x \in M \cap \dom (Pu)$.  Then $ax = sl$ for some $l \in
M_2$ and so $Pax = Psl$.  Using the property that $Pa - aP$ is zero on
$M_1$ and $Ps -sP$ is zero on $M_2$, we see that $aPx = sPl$ and we
deduce that $x \in \dom (uP)$.  Conversely if $x \in M \cap
\dom (uP)$, then $aPx = sPl$ for some $l \in L^2(G)$.  Using the
property that $Pa - aP$ is zero on $M_1$, we see that $ax =
P^{-1}sPl$.  But $aM_4 \cap P^{-1}sPL^2(G)
\subseteq sL^2(G)$, so $ax \in sL^2(G)$ and we deduce that $x \in
\dom(u)$.  Therefore
$M \cap \dom(Pu) = M\cap \dom(uP)$.  Finally for $x
\in M \cap \dom (u)$, we have $Pu x = y$ where $sP^{-1}y = ax$ and
$uPx = z$ where $sz = aP x$.  Thus $sP^{-1}y \in sM_2$ because
$aM_3 \cap sL^2(G) \subseteq sM_2$.  Since $\ker s = 0$ by Lemma
\ref{Ldense}, we see that $P^{-1}y \in M_2$ and hence $sy =
Pax$.  Also $Pax = aPx$ because $x \in M_1$.  Therefore $sy =
sz$ and since $\ker s = 0$ by Lemma \ref{Ldense}, we deduce
that $y = z$.  We conclude that $Pu = uP$ on $M$.

Conversely suppose there is a subspace $M$ of finite codimension in
$L^2(G)$ such that $Pu = uP$ on $M$.  Write $u = s^{-1}a = bt^{-1}$
where $a,b,s,t \in W(G)$ and let $N$ be a subspace of finite
codimension in $L^2(G)$ such that $tN \subseteq M$.  Then $tN \subset
\dom (u)$ and $sPu = suP$ on $M \cap \dom (u)$,
so $sPbn = aPtn$ for all $n \in N$ and we deduce that $sPb - aPt$ has
finite rank.  Furthermore $uP^{-1} = P^{-1}u$ on $P M$, a subspace of
finite codimension in $L^2(E) \oplus \mathbb {C}$, so by a similar
argument we see that $sP^{-1}b -aP^{-1}t$ also has finite rank.
Therefore $u \in R'(G)$ and we conclude from Theorem \ref{Tmain} that
$u \in D(G)$, as required.
\end{proof}

\providecommand{\bysame}{\leavevmode\hbox to3em{\hrulefill}\thinspace}


\begin{thebibliography}{1}

\bibitem{ber}
S.~K. Berberian, \emph{The maximal ring of quotients of a finite von {N}eumann
  algebra}, Rocky Mountain J. Math. \textbf{12} (1982), no.~1, 149--164.

\bibitem{berbook}
Sterling~K. Berberian, \emph{Baer {$*$}-rings}, Springer-Verlag, New York,
  1972, Die Grundlehren der mathematischen Wissenschaften, Band 195.

\bibitem{cohn}
P.~M. Cohn, \emph{Free rings and their relations}, second ed., Academic Press
  Inc. [Harcourt Brace Jovanovich Publishers], London, 1985.

\bibitem{cohn2}
\bysame, \emph{Skew fields}, Cambridge University Press, Cambridge, 1995,
  Theory of general division rings.

\bibitem{connes}
Alain Connes, \emph{Noncommutative geometry}, Academic Press Inc., San Diego,
  CA, 1994.

\bibitem{dieu2}
J.~Dieudonn{\'e}, \emph{Treatise on analysis. {V}ol. {I}{I}}, Academic Press
  [Harcourt Brace Jovanovich Publishers], New York, 1976, Enlarged and
  corrected printing, Translated by I. G. Macdonald, With a loose erratum, Pure
  and Applied Mathematics, 10-II.

\bibitem{duchreut}
G{\'e}rard Duchamp and Christophe Reutenauer, \emph{Un crit\`ere de
  rationalit\'e provenant de la g\'eom\'etrie non commutative}, Invent. Math.
  \textbf{128} (1997), no.~3, 613--622.

\bibitem{zero}
P.~A. Linnell, \emph{Zero divisors and group von {N}eumann algebras}, Pacific
  J. Math. \textbf{149} (1991), no.~2, 349--363.

\bibitem{div}
Peter~A. Linnell, \emph{Division rings and group von {N}eumann algebras}, Forum
  Math. \textbf{5} (1993), no.~6, 561--576.

\end{thebibliography}
\end{document}